\documentclass[12pt,twoside]{article}
\usepackage{graphicx}
\usepackage{graphics}
\usepackage{amsmath,amssymb,amsthm}

\theoremstyle{plain}
\newtheorem{thm}{Theorem}[section]

\newtheorem{lem}{Lemma}[section]
\newtheorem{cor}{Corollary}[section]
\theoremstyle{remark}

\newtheorem*{pf}{Proof}

\def\P{{\mathbb {P}}}                       
\def\E{{\mathbb {E}}}

\def\FD{{\mathcal F}}

\def\o{{\rm {o}}} 
\def\O{{\rm {O}}} 

\begin{document}
\centerline{\large{\bf {A population evolution model}}}
\smallskip
\centerline{\large{\bf {and its applications to random networks}}}

\bigskip
\centerline{ {\sc Istv\'an Fazekas},  {\sc Csaba Nosz\'aly} and {\sc Attila Perecs\'enyi} }

\medskip
\centerline{
University of Debrecen, Faculty of Informatics,} 
\centerline{
P.O. Box 400, 4002 Debrecen, Hungary,} 
\centerline{
e-mail: fazekas.istvan@inf.unideb.hu}

\bigskip
\begin{abstract}
A general population evolution model is considered. 
Any individual of the population is characterized by its score. 
Certain general conditions are assumed concerning the number of the individuals and their scores.
Asymptotic theorems are obtained for the number of individuals having some fixed score.
It is proved that the score distribution is scale free.
The result is applied to obtain the weight distributions of the cliques in a random graph evolution model.
\end{abstract}
\renewcommand{\thefootnote}{}
{
%
\section{Introduction}  \label{SectInt}
To describe real-life networks such as the WWW, social and biological networks, the preferential attachment model was introduced by Barab\'asi and Albert \cite{barabasi}. 
Then it was proved that the preferential attachment model leads to a scale-free random graph (for a rigorous mathematical proof see Bollob\'as et al. \cite{bollobas}). 
A random graph is called scale-free if it has a power law (asymptotic) degree distribution. 
Following the paper of Barab\'asi and Albert \cite{barabasi} several versions of the preferential attachment model were proposed. 
For the theory of random graphs one can consult with the monographs \cite{durrett}, \cite{janson}, \cite{hoff}.

In Ostroumova et al. \cite{ORS} a general graph evolution scheme was presented which covers lot of preferential attachment type models. 
They define a PA-class which covers the original preferential attachment model, the Holme-Kim model, the random Apollonian network, and the Buckley-Osthus model.
They  proved that the PA-class model leads to a scale-free graph.

In this paper we present a further generalization of the model by Ostroumova et al. \cite{ORS}. 
It is well-known that population evolution models and random graph evolution models are closely related.
So we consider the evolution of a population where any individual is characterized by its score.
We call it Model S.
During the evolution both the size of the population and the scores of the individuals can be increased. 
Let $X_n(s)$ denote the number of individuals having score $s$ at time $n$. 
First we describe the behaviour of the expectation  $\E X_n(s)$, see Theorem~\ref{ExpThm}.
Then we prove that the score distribution is scale-free, Theorem~\ref{limThm} and  Corollary~\ref{scalefreeCorr}. 
Our results generalize those of \cite{ORS}.
To obtain the above results we apply the methods presented in \cite{ORS}.

Then we apply our results to a random graph which is based on $N$-interactions.
The $N$-interactions model for $N=3$ was introduced in Backhausz and M\'ori \cite{BaMo1} (see also \cite{BaMo2}).
The general $N$-interactions model introduced and studied in  Fazekas and Porv\'azsnyik \cite{FIPB}, \cite{FIPB2}, \cite{FIPB3}.
The model incorporates the preferential attachment rule and the uniform choice of vertices.
In that model the vertices and the cliques possess certain weights.
In this paper we obtain that in the $N$-interactions model the weight distribution of the cliques is a power law, Theorem~\ref{CliqueThm}.
This theorem  generalizes the results of Fazekas et al. \cite{FNP} where the case $N=3$ was covered.

\section{Model S and its asymptotic behaviour} \label{SectModelS}
\setcounter{equation}{0}
We describe the evolution of a population. 
The evolution procedure defined below till equation \eqref{Exi} is called model S.
The evolution starts at time $0$ with maximum $t$ individuals. 
At each time $n=1,2,\dots$ maximum $t$ individuals are born. 
Any individual is characterized by its score. 
At birth the score is $u$ with high probability. 
More precisely the score of a new individual is at least $u$ and at time $n$ 
\begin{equation}  \label{score>u}
\P\left(\hbox{the score of a new individual}>u\right)=\O\left(\frac{1}{n}\right).
\end{equation}
The score of the individual $i$ at time $n$ is denoted by $S_n(i)$. 
Let $\FD_n$ denote the past of the population up to time $n$. 
The evolution of the score is described by the following equations
\begin{equation*}
\P\left(S_{n+1}(i)=S_n(i)+1|\FD_n\right)=a\frac{S_n(i)}{n}+b\frac{1}{n}+\O\left(\left(\frac{S_n(i)}{n}\right)^2\right),
\end{equation*}
\begin{equation}  \label{BasicEq}
\P\left(S_{n+1}(i)=S_n(i)|\FD_n\right)=1-a\frac{S_n(i)}{n}-b\frac{1}{n}+\O\left(\left(\frac{S_n(i)}{n}\right)^2\right),
\end{equation}
\begin{equation*}
\P\left(S_{n+1}(i)>S_n(i)+1|\FD_n\right)=\O\left(\left(\frac{S_n(i)}{n}\right)^2\right),
\end{equation*}
where $a$ and $b$ are fixed non-negative numbers. 
So at each step the score is increased by $1$ or $0$, the higher increasing is of low probability. 
Assume also that the total increase of the scores is at most $t$ at each step. 
Denote by $\xi_n$ the number of new individuals at time $n$. 
Assume that
\begin{equation}  \label{Exi}
\E\xi_n=m+\O\left(\frac{1}{n}\right)
\end{equation}
where $m>0$.

Let $X_n(s)$ denote the number of individuals having score $s$ at time $n$. 
Let $\Theta(x)$ denote a quantity with $|\Theta(x)|<x$. 
The first theorem shows that the expectation of the score distribution is scale-free.
\begin{thm}  \label{ExpThm}
Suppose that the conditions of Model S are satisfied and $a>0$. 
Then for any fixed $s=u,u+1,u+2,\dots$
\begin{equation}   \label{main}
\E X_n(s)=c(u,s)\left(n+\Theta\left(Ks^{2+\frac{1}{a}}\right)\right)
\end{equation}
for all $n$, where $K$ is a fixed finite constant,
\begin{equation}  \label{cus}
c(u,s)=\frac{\varGamma\left(s+\displaystyle{\frac{b}{a}}\right)\varGamma\left(u+\displaystyle{\frac{b+1}{a}}\right)}{a\varGamma\left(s+\displaystyle{\frac{b+a+1}{a}}\right)\varGamma\left(u+\displaystyle{\frac{b}{a}}\right)}m
\end{equation}
and $\varGamma$ denotes the $\varGamma$-function. 
Moreover
\begin{equation}  \label{cusGamma}
c(u,s)\sim\frac{m\varGamma\left(u+\displaystyle{\frac{b+1}{a}}\right)}{a\varGamma\left(u+\displaystyle{\frac{b}{a}}\right)}s^{-1-\frac{1}{a}} \ \ \hbox{ as } \ \ s\rightarrow\infty.
\end{equation}
\end{thm}
Now we consider the concentration of the scores around the expectation.
\begin{thm} \label{limThm}
Suppose that the conditions of Model S are satisfied and $a>0$.
Then there exists a constant $C$ such that for any fixed $s$ the following holds:
\begin{equation}
\P\left(|X_n(s)-\E X_n(s)|\geq s\sqrt{n}\log n\right) \le 2 n^{-\log n/C}.
\end{equation}
Furthermore let $\delta>0$. 
Then there exists a function $\varphi(n) = \o(1)$ such that with probability 1 we have: 
there exists only finitely many values of $n$ such that 
$|X_n(s)-\E X_n(s)|\geq\varphi(n)\E X_n(s)$
for any $s\leq n^{\frac{a-\delta}{4a+2}}$.
\end{thm}
Now we find that the score distribution is scale-free.
\begin{cor} \label{scalefreeCorr}
Suppose that the conditions of Model S are satisfied and $a>0$.
Then
\begin{equation*}
\lim_{n\rightarrow\infty}\frac{X_n(s)}{n}=c(u,s)
\end{equation*}
with probability $1$, where $c(u,s)$ satisfies  \eqref{cusGamma}.
\end{cor}
\section{Proofs of the main theorems and auxiliary results}  \label{SectProof}
\setcounter{equation}{0}
\begin{pf}[of Theorem \ref{ExpThm}]
We shall use the following notation.
\begin{equation}  \label{qn1s}
q_n^{(1)}(s)=\P\{S_{n+1}(i)=s+1|S_n(i)=s\}=a\frac{s}{n}+b\frac{1}{n}+\O\left(\frac{s^2}{n^2}\right),
\end{equation}
\begin{equation}  \label{qns}
q_n^{(j)}(s)=\P\{S_{n+1}(i)=s+j|S_n(i)=s\}=\O\left(\frac{s^2}{n^2}\right),\ \ \hbox{ }j\geq 2,
\end{equation}
\begin{equation}   \label{qn}
q_n=\P\left(\hbox{the score of a new individual}>u\right)=\O\left(\frac{1}{n}\right),
\end{equation}
\begin{equation}   \label{qSum}
q_n(s)=\sum_{j=1}^{\infty}q_n^{(j)}(s).
\end{equation}
The above relations are just the assumptions of Model S.

Part A of the proof.
We see that $q_n(s)$ is the probability that the score of an individual of score $s$ is increased at time $n$.
By \eqref{qn1s}, we obtain
\begin{equation}  \label{qconnect}
\frac{as+b+1}{as-a+b}q_n^{(1)}(s-1)=q_n^{(1)}(s)+\frac{1}{n}+\O\left(\frac{s^2}{n^2}\right).
\end{equation}
We have to calculate the conditional expectations
\begin{equation}  \label{EXn+1u|}
\E\{X_{n+1}(u)|X_n(u)\}=X_n(u)(1-q_n(u))+\xi_n(1-q_n),
\end{equation}
\begin{equation}   \label{EXn+1s}
\E\{X_{n+1}(s)|X_n(s),X_n(s-1),\dots,X_n(s-t)\}=
\end{equation}
\begin{equation*}
=X_n(s)(1-q_n(s))+X_n(s-1)q_n^{(1)}(s-1)+\sum_{j=2}^t X_n(s-j)q_n^{(j)}(s-j)+\O(q_n).
\end{equation*}

%
Consider the quantity $c(u,s)$ defined by \eqref{cus}.
Using the Stirling formula, we can see that \eqref{cusGamma} is true.
By the assumptions of our model,  $X_n(s) s\leq n t$.
This fact and \eqref{cusGamma} imply that
\begin{equation}   \label{Xbound}
X_n(s)=\O(n c(u,s)s^{\frac{1}{a}}).
\end{equation}
Therefore 
\begin{equation*}
\left| \E X_n(s)-n c(u,s)\right|  \le  \tilde{K} n c(u,s)s^{\frac{1}{a}} + n c(u,s) = n c(u,s)  \left( \tilde{K} s^{\frac{1}{a}} +1 \right)
\le c(u,s) s^{2+ \frac{1}{a}} 2 C_7 \left( \tilde{K} +1 \right)
\end{equation*}
holds for $n<2C_7s^2$.
So \eqref{main} is true with $K \ge \hat{K} = 2 C_7 \left( \tilde{K} +1 \right)$
if
\begin{equation}   \label{C7} 
n<2C_7s^2,
\end{equation}
 where $C_7$ will be defined later.

We shall use mathematical induction on $s$.

Part B of the proof.
First we study the case of the smallest value of $s$ that is when $s=u$. 
We see that the value of $c(u,s)$ if $s=u$ is
\begin{equation*}
c(u,u)=\frac{m}{au+b+1}.
\end{equation*}
Therefore to obtain \eqref{main} for $s=u$ we have to prove
\begin{equation}   \label{EXnu}
\E X_n(u)=\frac{nm}{au+b+1}+\Theta(C_1).
\end{equation}
For small values of $n$ we have already proved \eqref{EXnu}.
Assume now that \eqref{EXnu} is true for a certain $n$.
Now using \eqref{EXn+1u|} and \eqref{EXnu}, then applying \eqref{qn1s} and \eqref{qns}, we obtain
\begin{equation*}
\E X_{n+1}(u)=\E X_n(u)\left(1-q_n(u)\right)+\E\xi_n(1-q_n)=
\end{equation*}
\begin{equation*}
= \left(\frac{nm}{au+b+1}+\Theta(C_1)\right)\left(1-q_n(u)\right)+\E\xi_n(1-q_n)=
\end{equation*}
\begin{equation*}
=\frac{(n+1)m}{au+b+1}-\frac{m}{au+b+1}+\Theta(C_1)\left(1-q_n(u)\right)-q_n(u)\frac{nm}{au+b+1}+\E\xi_n(1-q_n)=
\end{equation*}
\begin{equation*}
=\frac{(n+1)m}{au+b+1}+\Theta(C_1)\left(1-q_n(u)\right)-\frac{m}{au+b+1}-\frac{\left(\frac{au+b}{n}+\O\left(\frac{u^2}{n^2}\right)\right)nm}{au+b+1}+\E\xi_n(1-q_n)=
\end{equation*}
\begin{equation*}
=\frac{(n+1)m}{au+b+1}+\Theta(C_1)\left(1-q_n(u)\right)-\frac{m}{au+b+1}-\frac{(au+b)m}{au+b+1}-\frac{\O\left(\frac{u^2}{n}\right)}{au+b+1}+\E\xi_n(1-q_n)=
\end{equation*}
\begin{equation*}
=\frac{(n+1)m}{au+b+1}+\Theta(C_1)\left(1-q_n(u)\right)-\frac{\O\left(\frac{u^2}{n}\right)}{au+b+1}+\O\left(\frac{1}{n}\right).
\end{equation*}
In the last step we applied that $\E\xi_n(1-q_n) = m+\O\left(\frac{1}{n}\right)$.
By the above inequality
\begin{equation}  \label{EXn+1u}
\E X_{n+1}(u)=\frac{(n+1)m}{au+b+1}+\Theta(C_1),
\end{equation}
if
\begin{equation*}
C_1q_n(u)\geq\frac{C_3}{n}\frac{1}{au+b+1}+\frac{C_2}{n}.
\end{equation*}
This inequality is equivalent to 
\begin{equation*}
C_1\left(\frac{au+b}{n}+\O\left(\frac{u^2}{n^2}\right)\right)>\frac{C_3}{n}\frac{1}{au+b+1}+\frac{C_2}{n}
\end{equation*}
and
\begin{equation*}
C_1\left(au+b\right)>\frac{C_1c_0u^2}{n}+\frac{C_3}{au+b+1}+C_2.
\end{equation*}
This last equality holds for large $n$ and $C_1$.

So we have obtained that the induction step \eqref{EXnu} $\Rightarrow$ \eqref{EXn+1u} is true if $n\ge n_1$ and $C_1$ is large enough.
Now choose $C_7$ in \eqref{C7} so that  
\begin{equation}  \label{C7+}
n_1 <2C_7 u^2.
\end{equation}
Therefore in the case of $s=u$ we obtained that relation \eqref{main} is true with $K\ge \max \left\{ \hat{K}, C_1\right\}$.

%
Part C of the proof.
Let $s>u$ and assume that \eqref{main} holds for all scores smaller than $s$. 
For fixed $s$, we apply induction on $n$.
At the beginning of the proof we proved that \eqref{main} is true for $n<2 C_7 s^2$.
So assume that \eqref{main} is satisfied for some $n$.
To avoid confusion we mention the following.
For large values of $s$ we shall find that the induction goes with appropriately large but fixed $K$.
However, for small values of $s$ the induction step is true for an appropriate increasing sequence $K_s$.
The largest of the above mentioned $K$-values will fit to our goal. 

By \eqref{EXn+1s} and using induction
\begin{equation*}
\E X_{n+1}(s)=\E X_n (s)\left(1-q_n(s)\right)+\E X_n (s-1) q^{(1)}_n(s-1)+\sum_{j=2}^t\E X_n(s-j)q_n^{(j)}(s-j)+\O(q_n)=
\end{equation*}
\begin{equation*}
=c(u,s)\left(n+\Theta(K_s s^{2+\frac{1}{a}})\right)\left(1-q_n(s)\right)+c(u,s-1)\left(n+\Theta(K_{s-1} (s-1)^{2+\frac{1}{a}})\right)q_n^{(1)}(s-1)+
\end{equation*}
\begin{equation*}
+\Theta\left(\frac{C_4c(u,s)s^{2+\frac{1}{a}}}{n}\right).
\end{equation*}
In the last step we applied \eqref{Xbound} and \eqref{cusGamma}.
Now using
\begin{equation*}
c(u,s) = \frac{as-a+b}{as+b+1} c(u,s-1),
\end{equation*}
we obtain
\begin{equation*}
\E X_{n+1}(s)=
\end{equation*}
\begin{equation*}
=c(u,s)(n+1)+c(u,s-1)nq^{(1)}(s-1)-c(u,s)nq_n(s)-c(u,s)+c(u,s)\Theta(K_ss^{2+\frac{1}{a}})\left(1-q_n(s)\right)+
\end{equation*}
\begin{equation*}
+c(u,s)\frac{as+b+1}{as-a+b}\Theta(K_{s-1}(s-1)^{2+\frac{1}{a}})q_n^{(1)}(s-1)+\Theta\left(\frac{C_4c(u,s)s^{2+\frac{1}{a}}}{n}\right)=
\end{equation*}
\begin{equation*}
=c(u,s)(n+1)+c(u,s)\Theta(K_ss^{2+\frac{1}{a}})\left(1-q_n(s)\right)+
\end{equation*}
\begin{equation*}
+ c(u,s)\frac{as+b+1}{as-a+b}\Theta(K_{s-1}(s-1)^{2+\frac{1}{a}})q_n^{(1)}(s-1) +\Theta\left(\frac{C_5c(u,s)s^{2+\frac{1}{a}}}{n}\right).
\end{equation*}
At the last step we used that
\begin{equation*}
c(u,s)\frac{as+b+1}{as-a+b}nq_n^{(1)}(s-1)-c(u,s)nq_n(s)=c(u,s)\left(1+\O\left(\frac{s^2}{n}\right)\right)
\end{equation*}
which is valid because of \eqref{qconnect}, \eqref{qns} and the assumption that the total increase of the scores is at most $t$.

To finish the induction we have to show that there exist a $K=K_s$ constant such that
\begin{equation}\label{Kineq}
K_ss^{2+\frac{1}{a}}q_n(s)\geq K_{s-1}\frac{as+b+1}{as-a+b}(s-1)^{2+\frac{1}{a}}q_n^{(1)}(s-1)+C_5\frac{s^{2+\frac{1}{a}}}{n}.
\end{equation}
Using Taylor's expansion, we obtain
\begin{equation*}
(s-1)^{2+\frac{1}{a}}  \le s^{2+\frac{1}{a}} - \left(2+\frac{1}{a}\right)s^{1+\frac{1}{a}} + C_6 s^{\frac{1}{a}}.
\end{equation*}
Therefore we have to show that
\begin{equation*}
Ks^{2+\frac{1}{a}}q_n(s)\geq K\frac{as+b+1}{as-a+b}\left(s^{2+\frac{1}{a}}-\left(2+\frac{1}{a}\right)s^{1+\frac{1}{a}}+C_6s^{\frac{1}{a}}\right)q_n^{(1)}(s-1)+C_5\frac{s^{2+\frac{1}{a}}}{n}.
\end{equation*}
Using \eqref{qconnect} and \eqref{qn1s}, we have to show that
\begin{equation*}
K\frac{s^{2+\frac{1}{a}}}{n}\left(2a+\frac{(b+1)(2a+1)}{as}+\O\left(\frac{s}{n}\right)\right)\geq Ks^{2+\frac{1}{a}}\O\left(\frac{s^2}{n^2}\right)+
\end{equation*}
\begin{equation*}
+K\frac{as+b+1}{as-a+b}C_6s^{\frac{1}{a}}\left(a\frac{s-1}{n}+b\frac{1}{n}+\O\left(\frac{s^2}{n^2}\right)\right)+C_5\frac{s^{2+\frac{1}{a}}}{n}.
\end{equation*}
This is equivalent to the following:
\begin{equation}  \label{K++}
K\frac{s^{2+\frac{1}{a}}}{n}\geq KC_7\frac{s^{4+\frac{1}{a}}}{n^2}+KC_8\frac{s^{1+\frac{1}{a}}}{n}+C_9\frac{s^{2+\frac{1}{a}}}{n} +
\end{equation}
\begin{equation*}
+ KC_{10}\frac{s^{3+\frac{1}{a}}}{n^2}+KC_{11}\frac{s^{2+\frac{1}{a}}}{n^2}+KC_{12}\frac{s^{\frac{1}{a}}}{n}.
\end{equation*}
Here $C_7, \dots , C_{12} $ are appropriate constants satisfying also $C_7 \ge \frac{n_1}{2u^2}$ (see \eqref{C7} and \eqref{C7+}).
So we have to show that 
\begin{equation*}  
K s^{2}\geq KC_7\frac{s^{4}}{n}+KC_8s+C_9 s^{2} + KC_{10}\frac{s^{3}}{n}+KC_{11}\frac{s^{2}}{n}+KC_{12}.
\end{equation*}
So for $n> 2 C_7 s^2$ we have to show that
\begin{equation*}  
K s^{2}\geq K\frac{s^{2}}{2}+KC_8s+ C_9 s^{2} + KC_{10}\frac{s}{2C_7}+KC_{11}\frac{1}{2C_7}+KC_{12}.
\end{equation*}
It is true if $K>C_{13}$ and $s>s_1$.
So \eqref{K++} holds, if $n\geq 2 C_7s^2$, $K>C_{13}$ and $s>s_1$.

We have already proved the desired result for $n <C_7s^2$.
Only the case of $s\le s_1$ is left.
However, as the magnitude of $q_n(s)$ is the same as that of $q^{(1)}_n(s-1)$, inequality \eqref{Kineq} is true for an appropriate increasing sequence $K_s$.
So part C of the proof is complete if we choose $K\ge \hbox{max}\{K_{s_1},C_{13}\}$.

Finally, we have to choose $K=\hbox{max}\{K_{s_1},C_{13}, \hat{K}, C_1\}$, where $\hat{K}$ and $C_1$ are from parts A and B of the proof. 
\qed
\end{pf}
\begin{lem} \label{lem1}
Let $\tau$ is a random variable and let $\eta_1$, $\eta_2$ be stochastic elements such that $\eta_1$ is $\eta_2$-measurable. Then
\begin{equation*}
|\E(\tau|\eta_2)-\E(\tau|\eta_1)|\leq \max_y\E(\tau|\eta_2=y)-\min_y\E(\tau|\eta_2=y)
\end{equation*}
almost surely.
\end{lem}
\begin{pf}
Let $g(y)= \E(\tau|\eta_2=y)$. 
Then
\begin{equation*}
\E(\tau|\eta_2)=g(\eta_2),
\end{equation*}
and therefore
\begin{equation*}
\E(\tau|\eta_1)=\E\{\E(\tau|\eta_2)|\eta_1\}=\E(g(\eta_2)|\eta_1).
\end{equation*}
So we have to prove that
\begin{equation*}
|g(\eta_2)-\E\{g(\eta_2)|\eta_1\}|\leq \max_y g(y)-\min_y g(y).
\end{equation*}
It is obvious that
\begin{equation*}
\min_y g(y)\leq g(\eta_2)\leq\max_y g(y)
\end{equation*}
and
\begin{equation*}
-\max_y g(y)\leq -\E\{g(\eta_2)|\eta_1\}\leq-\min_y g(y)
\end{equation*}
almost surely.
The result follows if we add the above two inequalities.
\qed
\end{pf}
\begin{lem} \label{lem2}
Let $\tau$ be random variable.
Let $\xi_t$, $t\ge 0$, be a stochastic process, $t_1,t_2$ be numbers with $0\le t_1 \le t_2$.
Let $\eta_1$, $\eta_2$ be stochastic elements defined by 
$\eta_1 = \{ \xi_t \, : \, 0 \le t \le t_1 \}$, $\eta_2 = \{ \xi_t \, : \, 0 \le t \le t_2 \}$.
That is $\eta_1$ and  $\eta_2$ are restrictions of the process $\xi_t$ up to time $t_1$ and $t_2$, respectively.
We see that $\eta_1$ is $\eta_2$-measurable.
If $y$ denotes an arbitrary realization of $\eta_2$ (that is $y= \eta_2(\omega)$ for some fixed $\omega \in \Omega$), 
then we denote by $\tilde{y}$ the restriction of $y$ up to time $t_1$.
Then
\begin{equation*}
|\E(\tau|\eta_2=y_2)-\E(\tau|\eta_1=y_1)|  \leq   \max_{\{y \,: \,\tilde{y}=y_1\}}\E(\tau|\eta_2=y) - \min_{\{y \, : \, \tilde{y}=y_1\}} \E(\tau|\eta_2=y)
\end{equation*}
almost surely, if $\tilde{y}_2=y_1$.
\end{lem}
\begin{pf}
Let $g(y)=\E(\tau|\eta_2=y)$. 
Then
\begin{equation*}
\E(\tau|\eta_2)=g(\eta_2).
\end{equation*}
As
\begin{equation*}
\E(\tau|\eta_1)=\E\{\E(\tau|\eta_2)|\eta_1\}=\E(g(\eta_2)|\eta_1),
\end{equation*}
so
\begin{equation*}
\E(\tau|\eta_1=y)=\E(g(\eta_2)|\eta_1=y).
\end{equation*}
Therefore we have to prove that, in case of $\tilde{y}_2=y_1$,
\begin{equation*}
|g(y_2)-\E(g(\eta_2)|\eta_1=y_1)|\leq \max_{\{y:\tilde{y}=y_1\}} g(y)-\min_{\{y:\tilde{y}=y_1\}} g(y).
\end{equation*}
It is obvious that, if $\tilde{y}_2=y_1$, then
\begin{equation*}
\min_{\{y:\tilde{y}=y_1\}} g(y)\leq g(y_2)\leq\max_{\{y:\tilde{y}=y_1\}} g(y)
\end{equation*}
and
\begin{equation*}
-\max_{\{y:\tilde{y}=y_1\}} g(y)\leq -\E\left(g(\eta_2)|\eta_1=y_1\right)\leq-\min_{\{y:\tilde{y}=y_1\}} g(y).
\end{equation*}
almost surely.
If we add the above two inequalities the result follows.
\qed
\end{pf}
\begin{pf}[of Theorem \ref{limThm}]
Let $Y_i=\E(X_n(s)|\FD_i)$, $i=0,1,\dots,n$, where $s$ and $n$ are fixed. 
It is easy to see that $(Y_i,\FD_i,i=0,1,\dots,n)$ is a martingale.  
Here $\FD_0$ is the trivial $\sigma$-algebra, so $Y_0=\E X_n(s)$. 
Furthermore $Y_n=X_n(s)$, because $X_n(s)$ is $\FD_n$-measurable. 
We will show that
\begin{equation} \label{Yineq}
|Y_{i+1}-Y_i|\leq Ms
\end{equation}
holds for some constant $M$.

Let $\eta_i$ denote the stochastic element describing the evolution of the population until the $i$th step. 
Then
\begin{equation*}
Y_i=\E(X_n(s)|\FD_i)=\E(X_n(s)|\eta_i)=g_i(\eta_i),
\end{equation*}
where $g_i(y)=\E(X_n(s)|\eta_i=y)$. 

During the following calculation when $\eta_i$ stands in a conditional expectation we shall assume that 
$\eta_i=y_i$ is fixed for any $i=0,1,\dots,n$.
Using Lemma \ref{lem2}, we obtain
\begin{equation}  \label{G}
|Y_{i+1}-Y_i|=|\E(X_n(s)|\eta_{i+1})-\E(X_n(s)|\eta_{i})|\leq g_{i+1}(\hat{y}_{i+1})-g_{i+1}(\bar{y}_{i+1}),
\end{equation}
where $\hat{y}_{i+1}=\arg\max g_{i+1}(y)$ and $\bar{y}_{i+1}=\arg\min g_{i+1}(y)$
and the maximum and the minimum are taken for $y$ having restriction $y_i$ up to time $i$. 
(Such $\hat{y}_{i+1}$ and $\bar{y}_{i+1}$ exist, because there are only finitely many orbits of population evolution until the $(i+1)$th step.)
Now we have to find bounds for
\begin{equation*}
\E(X_n(s)|\eta_{i+1}=\hat{y}_{i+1})-\E(X_n(s)|\eta_{i+1}=\bar{y}_{i+1})
\end{equation*}
in case of $0\leq i<n$.

Introduce notation
\begin{equation*}
\delta_l^{(i)}(s):=\E(X_l(s)|\eta_{i+1}=\hat{y}_{i+1})-\E(X_l(s)|\eta_{i+1}=\bar{y}_{i+1}),
\end{equation*}
where $i+1\leq l\leq n$.

In order to prove inequality \eqref{Yineq}, we need to show that $\delta_n^{(i)}(s)\leq Ms$ holds.

First we consider the case when $l\leq C_{14}s^2$. 
During one step of the evolution the total increase of the scores is at most $t$, so $s\delta_l^{(i)}(s)\leq tl$. 
So we obtain that
\begin{equation*}
\delta_l^{(i)}(s)\leq \frac{tl}{s}\leq\frac{tC_{14}s^2}{s}=C_{14}ts=Ms.
\end{equation*}

If $n=i+1$, then the only possible value of $l$ is $l=n=i+1$. 
In this case 
\begin{equation*}
\delta_{i+1}^{(i)}(s)=\E(X_{i+1}(s)|\eta_{i+1}=\hat{y}_{i+1})-\E(X_{i+1}(s)|\eta_{i+1}=\bar{y}_{i+1}).
\end{equation*}
However, here $\hat{y}_{i+1}$ and $\bar{y}_{i+1}$ are the same, except at the last time step, where they can be different. 
In one step the number of individuals of score $s$ can be changed at most by $t$. 
So we obtain
\begin{equation*}
\delta_{i+1}^{(i)}(s)\leq t\leq Ms.
\end{equation*}

Now we consider the case when $i\leq l\leq n-1$ and $l > C_{14}s^2$. 
Taking the conditional expectation with respect to  $\eta_i=y_i$ in \eqref{EXn+1u|} and \eqref{EXn+1s}, we obtain
\begin{equation}  \label{EXn+1ueta}
\E\{X_{l+1}(u)|\eta_i=y_i\}=\E\{X_l(u)|\eta_i=y_i\}(1-q_l(u))+\E\{\xi_l|\eta_i=y_i\},
\end{equation}
\begin{equation}   \label{EXn+1seta}
\E\{X_{l+1}(s)|\eta_i=y_i\}=
\end{equation}
\begin{equation*}
=\E\{X_l(s)|\eta_i=y_i\}(1-q_l(s))+\E\{X_l(s-1)|\eta_i=y_i\}q_l^{(1)}(s-1)+
\end{equation*}
\begin{equation*}
+\sum_{j=2}^t \E\{X_l(s-j)|\eta_i=y_i\}q_l^{(j)}(s-j)+\O(q_l)
\end{equation*}
for $s\ge u+1$.
Applying \eqref{G} for $l$ instead of $n$, relations  \eqref{qns}, \eqref{EXn+1ueta} and  \eqref{EXn+1seta} imply
\begin{equation} \label{deltaineq}
\delta_{l+1}^{(i)}(s)\leq \delta_{l}^{(i)}(s)(1-q_l(s))+\delta_{l}^{(i)}(s-1)q_l^{(1)}(s-1)+
\end{equation}
\begin{equation*}
+ \O\left(\frac{1}{l}\right)+\sum_{j=2}^t \delta_{l}^{(i)}(s-j)\O\left(\left(\frac{s-j}{l}\right)^2\right).
\end{equation*}

Now we use induction on $l$. 
We have already seen that if $l$ is small, then $\delta_{l}^{(i)}(s)<Ms$ holds. 
Using  \eqref{qn1s}-\eqref{qSum}, from \eqref{deltaineq} we get
\begin{equation*}
\delta_{l+1}^{(i)}(s)\leq Ms(1-q_l(s))+M(s-1)q_l^{(1)}(s-1)+\O\left(\frac{1}{l}\right)+\sum_{j=2}^t M(s-j)\O \left(\left(\frac{s-j}{l}\right)^2\right)=
\end{equation*}
\begin{equation*}
= Ms-Ms\left(\frac{as}{l}+\frac{b}{l}+\O\left(\frac{s^2}{l^2}\right)\right)+M(s-1)\left(\frac{a(s-1)}{l}+\frac{b}{l}+\O\left(\frac{s^2}{l^2}\right)\right)+
\end{equation*}
\begin{equation*}
+ \O\left(\frac{1}{l}\right)
+M\left(ts-\frac{t(t+1)}{2}\right)\O\left(\frac{s^2}{l^2}\right)=
\end{equation*}
\begin{equation*}
= Ms-M2s\frac{a}{l}+M\frac{a}{l}+M\frac{s^3}{l^2}C_{15}+O\left(\frac{1}{l}\right)+Mts\O\left(\frac{s^2}{l^2}\right)  \le Ms
\end{equation*}
if $l> C_{14} s^2$ with $C_{14} > C_{15}/a$ and with $M$ large enough.
So \eqref{Yineq} is proved.

As $Y_n = X_n(s)$, $Y_0 = \E X_n(s)$ and, by  \eqref{Yineq}, $|Y_{i+1} -Y_i| \le Ms$, so from the Azuma-Hoeffding inequality (see \cite{HO}, \cite{Azuma}) it follows that
\begin{equation} \label{AzHo}
\P(|X_n(s)-\E X_n(s)|\geq s\sqrt{n}\log n) = \P(|Y_n-Y_0|\geq s\sqrt{n}\log n)\leq 
\end{equation}
\begin{equation*}
\leq2\exp\{-\frac{s^2n\log^2n}{2nM^2s^2}\}=2(e^{\log n})^{-\log n/2M^2}=2n^{-\log n/2M^2}. 
\end{equation*}
Let $\delta>0$, and $s\leq n^{\frac{a-\delta}{4a+2}}$. Then
\begin{equation*}
\P\left\{|X_n(s)-\E X_n(s)|\geq\varphi(n)\E X_n(s)\right\} \leq
\end{equation*}
\begin{equation*}
\P\left\{|X_n(s)-\E X_n(s)|\geq\varphi(n)Cs^{-1-\frac{1}{a}}\left(n+\O\left(s^{2+\frac{1}{a}}\right)\right)\right\}.
\end{equation*}
We see that relation
\begin{equation} \label{phi}
\varphi(n)Cs^{-1-\frac{1}{a}}\left(n+\O\left(s^{2+\frac{1}{a}}\right)\right)>s\sqrt{n}\log n
\end{equation}
is equivalent to
\begin{equation*}
\varphi(n)C\left(s^{-2-\frac{1}{a}}n+\O(1)\right)>\sqrt{n}\log n.
\end{equation*}
Now in case of $s=n^{\frac{a-\delta}{4a+2}}$
\begin{equation*}
\varphi(n) C\left(s^{-2-\frac{1}{a}}n+\O(1)\right)  = \varphi(n)C\left(n^{\frac{a-\delta}{4a+2}\left(-\frac{2a+1}{a}\right)}n+\O(1)\right)=
\end{equation*}
\begin{equation*}
= \varphi(n)C\left(n^{-\frac{a-\delta}{2a}}n+\O(1)\right)
=\varphi(n)C\left(n^{1-\frac{a-\delta}{2a}}+\O(1)\right)=
\end{equation*}
\begin{equation*}
=\varphi(n)C\left(n^{\frac{a+\delta}{2a}}+\O(1)\right)=
\varphi(n)C\left(n^{\frac{1}{2}+\frac{\delta}{2a}}+\O(1)\right)>\sqrt{n}\log n
\end{equation*}
if $n$ is large and $\varphi(n)$ is a certain function with $\varphi(n)= \o(1)$.
So \eqref{phi} is valid for $s= n^{\frac{a-\delta}{4a+2}}$ therefore it is valid for $s\leq n^{\frac{a-\delta}{4a+2}}$.
Therefore the above calculation shows that for $s\leq n^{\frac{a-\delta}{4a+2}}$
\begin{equation} \label{AzHo+}
\P\left\{|X_n(s)-\E X_n(s)|\geq\varphi(n)\E X_n(s)\right\}  \leq  \P\left\{ |X_n(s)-\E X_n(s)|\geq s\sqrt{n}\log n\right\}. 
\end{equation}
So, by \eqref{AzHo} and by the Borel-Cantelli lemma with probability 1 we have: 
there exists only finitely many $n$ such that 
$|X_n(s)-\E X_n(s)|\geq\varphi(n)\E X_n(s)$
for any $s\leq n^{\frac{a-\delta}{4a+2}}$.
\qed
\end{pf}

\begin{pf}[of Corollary \ref{scalefreeCorr}]
From \eqref{AzHo} it follows that with probability $1$ the relation $\{|X_n(s)-\E X_n(s)|\geq s\sqrt{n}\log n\}$ holds only for finite number of $n$'s.
 So we obtain
\begin{equation*}
\frac{X_n(s)}{n}-\frac{\E X_n(s)}{n}\rightarrow 0 
\end{equation*}
with probability $1$.
From  Theorem \ref{ExpThm} we know that
\begin{equation*}
\frac{\E X_n(s)}{n}=c(u,s)\left(1+\frac{1}{n}\Theta\left(K s^{2+\frac{1}{a}}\right)\right)\rightarrow c(u,s)\ \ \hbox{ as } \ \ n\rightarrow\infty.
\end{equation*}
Therefore
\begin{equation*}
\lim_{n\rightarrow\infty}\frac{X_n(s)}{n}=c(u,s)
\end{equation*}
with probability $1$, and
\begin{equation*}
c(u,s)\sim\frac{m\varGamma\left(u+\displaystyle{\frac{b+1}{a}}\right)}{a\varGamma\left(u+\displaystyle{\frac{b}{a}}\right)}s^{-1-\frac{1}{a}} \ \ \hbox{ as } \ \ s\rightarrow\infty.
\end{equation*}
\qed
\end{pf}

\section{Application for $M$-cliques of the $N$-interactions model}  \label{SectAppl}
\setcounter{equation}{0}
In this section we apply the theorems of Section \ref{SectModelS} to the weights of the $M$-cliques of the $N$-interactions random graph model (see \cite{FIPB}, \cite{FIPB2}, \cite{FIPB3}). 
First we recall that, by the usual definition, a complete graph with $M$ vertices is called an $M$-clique. 
However, we have to emphasize that in our model only those complete graphs will be considered to be cliques which take their origin in interaction of $N$ vertices.
So we will see that during one step at most one new $N$-clique is constructed.

Let $N\geq 3$ and let $0<p\leq 1$, $0\leq r\leq 1$ and $0\leq q\leq 1$ be fixed numbers (the parameters of the $N$-interactions model). 
The evolution of our random graph starts at time $0$ with an $N$-clique. 
During the evolution any $N$-clique and any of its $M$-subclique $(0<M<N)$ has an initial weight $1$. 
(The weight of a non-existing clique is considered to be $0$.) 
In every time step $N$ vertices interact each other.
The interaction means that we draw all non-existing edges between them, so we construct an $N$-clique. 
Moreover, its weight and the weights of all of its subcliques are increased by $1$.

The details of the evolution are the following.
At each time step we have two possibilities. 
On the one hand, with probability $p$, a new vertex is added to the graph and it interacts with $N-1$ old vertices (so they form a new $N$-clique). 
On the other hand, with probability $(1-p)$, we do not add any new vertex but $N$ old vertices interact with each other (so they form an $N$-clique).

When a new vertex is added there are again two options. 
With probability $r$ we choose an $(N-1)$-clique according to the weights of the $(N-1)$-cliques, using the preferential attachment rule, while with probability $(1-r)$, we choose $N-1$ old vertices uniformly. 
(Here, the preferential attachment rule means that the $(N-1)$-clique $c$ of weight $w_c$ is chosen with probability $w_c/ \sum_k w_k$, where $\sum_k w_k$ is the total weight of all of the $(N-1)$-cliques.
Moreover, uniform choice means that $(N-1)$-cliques have equal chances among the $(N-1)$-cliques.)
Then the  $N-1$ old vertices chosen and the new vertex interact. 
It means that they form a new $N$-clique.
The other case, when we do not add any new vertex, there are again two possibilities. 
With probability $q$ we choose an $N$-clique according to the weights of the $N$-cliques, using the preferential attachment rule, while with probability $(1-q)$ we choose $N$ old vertices uniformly. 
The $N$ old vertices chosen interact, so they form an $N$-clique.

We stress that in any step and in any case the weight of the $N$-clique constructed in that step and the weights of its subcliques are increased by 1.

The following theorem shows that the asymptotic weight distribution of the $M$-cliques is power law.
The result was presented in \cite{FNP} for the case of $M=2$, $N=3$ and also for the case of $M=N$ for arbitrary $N>2$.
\begin{thm} \label{CliqueThm}
Let $N\ge 3$ be fixed and let $M$ be fixed with $2\le M\le N$ and denote by $X_M(n,w)$ the number of $M$-cliques having weight $w$ after $n$ steps.
If $p>0$ and either $r>0$ or $(1-p)q>0$, then
\begin{equation}
\dfrac{X_M \left( n,w \right)}{n} \rightarrow x_{M,w}
\end{equation}
almost surely, as $n \rightarrow \infty $, where $x_{M,w}$, $w=1,2, \dots$\,, are numbers satisfying
\begin{equation} \label{x_{M,w}_asz.}
x_{M,w} \sim \frac{\mu}{a} \varGamma \left( 1 + \frac{1}{a} \right) w^{- \left( 1 + \frac{1}{a} \right)}
\end{equation}
as $w \to \infty$, with
\begin{equation} \label{a}
a = p r \frac{N-M}{N}+ \left(1-p\right) q,
\end{equation}
and
\begin{equation} \label{mu}
\mu= p \binom{N-1}{M-1} + p(1-r) \binom{N-1}{M} + (1-p)(1-q) \binom{N}{M}.
\end{equation}
\end{thm}
\begin{pf}
Let $V_n$ denote the number of vertices after the $n$th step.
By the Marcinkiewicz strong law of large numbers, we have
\begin{equation}  \label{Markinkiewicz}
V_n = pn + \o \left( n^{1/2 + \varepsilon} \right)
\end{equation}
almost surely, for any $\varepsilon >0$.
We know that during one step the weight of an $M$-clique can be increased either by $1$ or by $0$.
We see that condition \eqref{score>u} is true with $u=1$.
Let $W_n(i)$ denote the weight of the $M$-clique $i$ at step $n$ ($2\le M\le N$).
Then for $w\ge 1$ we have
\begin{equation}  \label{Wn+1}
\P\left\{ W_{n+1}(i)=  w+1 \, | \, W_n(i)=w \right\} =
\end{equation}
\begin{equation*}
= pr \frac{(N-M)w}{Nn}  + p(1-r) \frac{\binom{V_n-M}{N-1-M} }{\binom{V_n}{N-1}} + (1-p) q \frac{w}{n} + (1-p)(1-q) \frac{\binom{V_n-M}{N-M} }{\binom{V_n}{N}} =
\end{equation*}
\begin{equation*}
= \left( p r \frac{N-M}{N}+ \left(1-p\right) q \right) \frac{w}{n} + \O\left(\frac{1}{n^M} \right) = a \frac{w}{n} + \O\left(\frac{1}{n^2} \right)
\end{equation*}
because of relation \eqref{Markinkiewicz} and condition $M\ge 2$.
So  we see that conditions \eqref{BasicEq} are true with $a$ from formula \eqref{a} and $b=0$.

Now let $\xi_n$ denote the number of the new $M$-cliques in step $n$.
Then 
\begin{equation*}
\E\left\{ \xi_{n+1} \, | \, \FD_n \right\} = pr \binom{N-1}{M-1} +
\end{equation*}
\begin{equation*}
+ p(1-r) \left[ \binom{N-1}{M-1} + \frac{\binom{N-1}{M} \binom{V_n}{N-1} -E_n \binom{V_n -M}{N-1-M} }{ \binom{V_n}{N-1}} \right]
+(1-p)(1-q) \left[ \binom{N}{M} - \frac{E_n \binom{V_n -M}{N-M} }{ \binom{V_n}{N}} \right],
\end{equation*}
where $E_n$ denotes the number of $M$-cliques at step $n$.
Therefore, applying $E_n \le n c$ and $M\ge 2$,
\begin{equation*}
\E\xi_{n+1}=\mu +\O\left(\frac{1}{n}\right),
\end{equation*}
where $\mu$ is defined in \eqref{mu}.
So condition \eqref{Exi} is satisfied.
Therefore Corollary \ref{scalefreeCorr} implies the result.
\qed
\end{pf}

\end{document}